\documentclass{amsart}

\usepackage{amsmath,amssymb,amscd,amsfonts}
\usepackage[all]{xy}

\newtheorem{theorem}{Theorem}
\newcommand{\bt}{\begin{theorem}}
\newcommand{\et}{\end{theorem}}
\newtheorem{lemma}{Lemma}
\newcommand{\bl}{\begin{lemma}}
\newcommand{\el}{\end{lemma}}
\newtheorem{corollary}{Corollary}
\newcommand{\bc}{\begin{corollary}}
\newcommand{\ec}{\end{corollary}}
\newcommand{\beq}{\begin{equation}}
\newcommand{\eeq}{\end{equation}}
\newcommand{\benum}{\begin{enumerate}}
\newcommand{\eenum}{\end{enumerate}}
\newcommand{\N}{\ensuremath{ \mathbf N }}
\newcommand{\Z}{\ensuremath{\mathbf Z}}
\newcommand{\Q}{\ensuremath{\mathbf Q}}
\newcommand{\R}{\ensuremath{\mathbf R}}

\newcommand{\mbF}{\ensuremath{ \mathbf F}}
\newcommand{\PP}{\ensuremath{ \mathbf P}}

\newcommand{\mcf}{\ensuremath{ \mathcal F}}

\newcommand{\mbx}{\ensuremath{ \mathbf x}}

\DeclareMathOperator{\card}{\text{card}}

\newcommand{\bmat}{\left(\begin{matrix}}
\newcommand{\emat}{\end{matrix}\right)}

\DeclareMathOperator{\qqand}{\qquad\text{and}\qquad}

\DeclareMathOperator{\vectorsmallxn}{\left( \begin{smallmatrix} x_1 \\ \vdots \\ x_n \end{smallmatrix}\right)}

\title{Sidon sets for linear forms}

\author{Melvyn B. Nathanson}
\address{Lehman College (CUNY), Bronx, New York 10468}
\email{melvyn.nathanson@lehman.cuny.edu}

\subjclass[2010]{11B13, 11B34, 11B75, 11P99}

\keywords{Sidon set,  sumset, sum of dilates, distinct subset sums, representation functions.}

\thanks{Supported in part by a grant from the PSC-CUNY Research Award Program.}

\date{\today}

\begin{document}

\begin{abstract}
Let $\varphi(x_1,\ldots, x_h) =  c_1 x_1 + \cdots + c_h x_h $ be a linear form 
with coefficients in a field \mbF, and let $V$ be a vector space over \mbF.  
A nonempty subset $A$  of $V$ is a 
\emph{$\varphi$-Sidon set} if 
$\varphi(a_1,\ldots, a_h) = \varphi(a'_1,\ldots, a'_h) 
$ implies $(a_1,\ldots, a_h) = (a'_1,\ldots, a'_h)$ 
for all $h$-tuples $(a_1,\ldots, a_h) \in A^h$ and $ (a'_1,\ldots, a'_h) \in A^h$.   
There exist infinite Sidon sets for the linear form $\varphi$ if and only 
if the set of coefficients of $\varphi$ has distinct subset sums.  
In a normed vector space with $\varphi$-Sidon sets, 
every infinite sequence of vectors is 
asymptotic to a $\varphi$-Sidon set of vectors.
Results on $p$-adic perturbations of $\varphi$-Sidon sets of integers and bounds on the growth 
of $\varphi$-Sidon sets of  integers are also obtained.  
\end{abstract}

\maketitle

\section{Linear forms with property $N$}               \label{LinearPerturb:section:forms}
Let \mbF\ be a field and let $h$ be a positive integer.  
We consider linear forms  
\beq                                                      \label{LinearPerturb:phi}
\varphi(x_1,\ldots, x_h) =  c_1 x_1 + \cdots + c_h x_h  
\eeq
where $c_i \in \mbF$ for all $i \in \{1,\ldots, h\}  $.  

Let $V$ be a vector space over the field \mbF.  
For every subset nonempty $A$ of $V$, let 
\[
A^h = \left\{ (a_1,\ldots, a_h): a_i \in A \text{ for all } i \in \{1,\ldots, h\}  \right\}
\]
be the set of all $h$-tuples of elements of $A$.  
For $c \in \mbF$, the \emph{$c$-dilate} of $A$ is the set 
\[
c\ast A = \{ca:a\in A\}.
\]
The \emph{$\varphi$-image of $A$} is the set 
\begin{align*}
\varphi(A) & = \left\{ \varphi(a_{1},\ldots, a_{h}): (a_1,\ldots, a_h) \in A^h \right\} \\ 
& = \left\{ c_1 a_{1} + \cdots +  c_h a_{h}: (a_1,\ldots, a_h) \in A^h \right\} \\
& = c_1\ast A + \cdots + c_h \ast A. 
\end{align*}
Thus, $\varphi(A)$ is a sum of dilates.  
We define $\varphi(\emptyset) = \{0\}$.

A nonempty subset $A$  of $V$ is a 
\emph{Sidon set for the linear form $\varphi$} or, simply, a 
\emph{$\varphi$-Sidon set} 
if it satisfies the following property:  
For all $h$-tuples $(a_1,\ldots, a_h) \in A^h$ and $ (a'_1,\ldots, a'_h) \in A^h$,   
if 
\[
\varphi(a_1,\ldots, a_h) = \varphi(a'_1,\ldots, a'_h) 
\]
then $(a_1,\ldots, a_h) = (a'_1,\ldots, a'_h)$, 
that is, $a_i = a'_i$ for all $i \in \{1,\ldots, h\} $.  
Thus, $A$ is a $\varphi$-Sidon set if the linear form $\varphi$ is one-to-one on $A^h$.  

Two cases of special interest are $V = \mbF$ with $\varphi$-Sidon sets contained in \mbF, 
and $V = \mbF = \Q$ with $\varphi$-Sidon sets of positive integers.

For the linear form $\varphi = \sum_{i=1}^h c_i x_i$, 
every set with one element is a $\varphi$-Sidon set.
There is a simple obstruction to the existence of $\varphi$-Sidon sets 
 with more than one element.  For every nonempty subset $I$ of $\{1,\ldots, h\}$, 
 define the \emph{subset sum} 
 \beq                                                            \label{LinearPerturb:sI}
 s(I) = \sum_{i \in I} c_i. 
\eeq
Let  $s({\emptyset}) = 0$.  
Suppose there exist disjoint subsets $I_1$ and $I_2$ of $\{1,\ldots, h\}$ 
with $I_1$ and $I_2$ not both empty such that 
\beq                                                      \label{LinearPerturb:obstruction}
s({I_1}) = \sum_{i \in I_1} c_{i} = \sum_{i \in I_2} c_{i} = s({I_2}).   
\eeq
Let $I_3 = \{1,\ldots, h\}  \setminus (I_1 \cup I_2)$.  
Let $A$ be a subset of $V$ with  $|A| \geq 2$.  Choose vectors $u,v,w  \in A$ with $u \neq v$, 
and define
\[
a_i = \begin{cases}
u & \text{ if $i \in I_1$}\\
v & \text{ if $i \in I_2$} \\
w & \text{ if $i \in I_3$}
\end{cases}
\]
and 
\[
a'_i = \begin{cases}
v & \text{ if $i \in I_1$ } \\
u & \text{ if $i \in I_2 $}  \\
w & \text{ if $i \in I_3$.}
\end{cases}
\]
We have 
\[
(a_1,\ldots, a_h) \neq (a'_1,\ldots, a'_h)
\]
because $I_1 \cup I_2 \neq \emptyset$ and $a_i \neq a'_i$ for all $i \in I_1 \cup I_2$.  

The sets $I_1$, $I_2$, $I_3$ are pairwise disjoint.   
Condition~\eqref{LinearPerturb:obstruction} implies 
\begin{align*}
\varphi(a_1,\ldots, a_h) 
& = \sum_{i\in I_1} c_i a_i + \sum_{i \in I_2} c_i a_i + \sum_{i\in I_3} c_i a_i \\
& = \left(  \sum_{i\in I_1} c_i \right) u  + \left( \sum_{i\in I_2} c_i  \right) v +  \left( \sum_{i\in I_3} c_i \right) w \\
& = \left( \sum_{i\in I_2} c_i \right) u  +  \left( \sum_{i\in I_1} c_i  \right) v+  \left( \sum_{i\in I_3} c_i \right) w \\
& =   \sum_{i\in I_1} c_i a'_i +   \sum_{i\in I_2} c_i  a'_i+  \sum_{i\in I_3} c_i a'_i \\
& = \varphi(a'_1,\ldots, a'_h)  
\end{align*} 
and so $A$ is not a $\varphi$-Sidon set.  

We say that the linear form~\eqref{LinearPerturb:phi} has \emph{property $N$}  
if there do \emph{not} exist disjoint subsets $I_1$ and $I_2$ of  $\{1,\ldots, h\}$ 
that satisfy condition~\eqref{LinearPerturb:obstruction} 
with $I_1$ and $I_2$ not both empty.  
If the linear form $\varphi = \sum_{i=1}^h c_ix_i$ has property $N$, 
then 
\[
\sum_{i\in I_1} c_i = s(I_1)  \neq s({\emptyset} )= 0
\]
for every nonempty subset $I_1$ of $\{1,\ldots, h\}$.  
In particular, choosing $I_1 = \{i\}$ shows that $c_i \neq 0$ for all $i \in \{1,\ldots, h\}$.

For example, if  $h \geq 1$ and  $c_i = 2^{i-1}$ for all $i \in \{1,\ldots, h\}  $,  
then the linear form 
\[
\varphi = \sum_{i=1}^h c_i x_i = x_1 + 2x_2 + 4x_3 + \cdots + 2^{h-1} x_h 
\]
has property $N$.  

If $h \geq 2$ and $c_i = 1$ for all $i \in \{1,\ldots, h\}$, then  the linear form   
\[
\psi = \sum_{i=1}^h c_i x_i = x_1 + x_2 + x_3 + \cdots + x_h  
\]  
does not have property $N$ because the nonempty disjoint sets $I_1 = \{1\}$ and $I_2 = \{2\}$ 
satisfy 
\[
\sum_{i \in I_1} c_{i}  = c_1 = 1 = c_2  = \sum_{i \in I_2} c_{i}.  
\]

In Section~\ref{LinearPerturb:section:contractions} we prove that, 
for every infinite vector space $V$, 
there exist infinite $\varphi$-Sidon sets for the linear form $\varphi$ 
if and only if $\varphi$ has property $N$.

For related work on additive number theory for linear forms, 
see ~Bukh\cite{bukh08} and Nathanson~\cite{nath2007-118,nath2007-120,nath2008-126,nath2009-133,nath2018-174,nath2007-119}.

Let $\varphi(x_1,\ldots, x_h) =  c_1 x_1 + \cdots + c_h x_h $, 
where $c_i \in \mbF$ for $i \in \{1,2,\ldots,h\}$.  
Let $J_1$ and $J_2$ be distinct subsets of $\{1,2,\ldots,h\}$ such that 
$\sum_{i\in J_1} c_i =  \sum_{i\in J_2} c_i$ and let $J = J_1 \cap J_2$.  
The  sets $I_1 = J_1\setminus J$ and  $I_2 = J_2 \setminus J$ are distinct and disjoint subsets 
of $\{1,2,\ldots,h\}$.  Moreover, $\sum_{i\in I_1} c_i =  \sum_{i\in I_2} c_i$. 
It follows that the linear form $\varphi$ has property $N$ if and only if the set 
$\{c_1,\ldots, c_h\}$ has distinct subset sums.  

Let $g(n)$ be the size of the largest subset of $\{1,2,\ldots,n\}$ that has 
distinct subset sums.  A famous unsolved problem of Paul Erd\H os and Leo Moser 
asks if 
\[
g(n) = \frac{\log n}{\log 2
} + O(1).
\]
See Erd\H os~\cite[pp. 136--137]{erdo56}, Guy~\cite[Section C8]{guy04}, 
and Dubroff, Fox, and Xu~\cite{dubr-fox-xu21}.

\section{Classical Sidon sets}                 \label{LinearPerturb:section:ClassicalSidon}

The idea of a Sidon set for a linear form derives from the classical definition of a Sidon set of integers.  
In additive number theory,  a  \emph{Sidon set} (also called a $B_2$-set) is a set $A$ of positive integers 
such that, if $a_1,a_2,a'_1,a'_2 \in A$ and 
\[
a_1 + a_2= a'_1 + a'_2
\]
then $\{a_1,a_2 \} = \{a'_1, a'_2\}$.  
More generally, let $G$ be an additive abelian group or semigroup, 
and let $A$ be a subset of $G$.
For $h \geq 2$, the \emph{$h$-fold sumset} of $A$ is the set $hA$ 
of all sums of $h$ not necessarily distinct elements of $A$.
A nonempty set $A$ is an \emph{$h$-Sidon set} (or a $B_h$-set) 
if every element of the sumset $hA$ has an 
essentially unique representation as the sum of $h$ elements of $A$, in the following sense: 
If  $\{a_i:i\in I\}$ is a set of pairwise distinct elements of $A$  
and if $\{u_i:i\in I\}$ and $\{v_i:i\in I\}$ are sets of nonnegative integers such that 
\[
h = \sum_{i\in I} u_i  = \sum_{i\in I} v_i  
\]
and
\[
\sum_{i\in I} u_ia_i = \sum_{i\in I} v_ia_i 
\]
then $u_i = v_i$ for all $i \in I$.

The sumset $hA$ is associated with the linear form 
\[
\psi = \psi(x_1,\ldots, x_h) = x_1 + \cdots + x_h  
\]
and  
\[
hA = \psi(A) 
 = \left\{ a_{1} + \cdots +  a_{h}: a_i \in A \text{ for all } i \in \{1,\ldots, h \} \right\}. 
\] 
The linear form $\psi$ does not satisfy condition $N$, and there  exists 
no $\psi$-Sidon set $A$ with $\card(A) \geq 2$.  

The literature on classical Sidon sets is huge.
Two surveys of results on classical Sidon sets are Halberstam and Roth~\cite{halb-roth66} 
and O'Bryant~\cite{obry04}. 
For recent work, see~\cite{cill-serr-wotz20,kiss-sand21,koha-lee-more-rol18,kova-mlad-tan17,liu-pach19,pach19,pach-sand19,scho-shkr19,xu18}.

\section{Contractions of linear forms}                                     \label{LinearPerturb:section:contractions}
Let \mbF\ be a field 
and let $\varphi = \sum_{i=1}^h c_ix_i$ be a linear form 
in $h$ variables with coefficients $c_i \in \mbF$.   
Associated to every subset $J$  of  $\{1,\ldots, h\}$ 
is the linear form in $\card(J)$ variables 
\[
\varphi_J = \sum_{j \in J} c_jx_j.
\]
We have $\varphi_{\emptyset} = 0$ and $\varphi_J = \varphi$ if  $J = \{1,\ldots, h\}$.   
The linear form $\varphi_J$ is called a \emph{contraction} of the linear form $\varphi$.

Let $V$ be a vector space over the field \mbF.  
For every nonempty subset  $A$ of $V$, let 
\[
\varphi_J(A) =  \left\{ \sum_{j \in J} c_j a_j : a_j \in A \text{ for all } j \in J \right\}. 
\]
If $A$ is a $\varphi$-Sidon set, then $A$ is a  $\varphi_J$-Sidon set 
for every nonempty subset $J$ of  $\{1,\ldots, h\}$.

For every subset $X$ of $V$ and vector $v \in V$, 
 the \emph{translate} of $X$ by $v$ is the set 
\[
X+v = \{ x+v:x\in X\}.  
\]
For every subset of $J$ of $\{1,\ldots, h\}$, let $J^c = \{1,\ldots, h\} \setminus J$ 
be the complement of $J$ in $\{1,\ldots, h\}$.   
For every subset $A$ of $V$ and $b \in V \setminus A$, we define 
\beq                                                                              \label{LinearPerturb:PhiAb}
\Phi_J(A,b) = \varphi_J(A) + \left( \sum_{j \in J^c} c_j \right) b = \varphi_J(A) +  s(J^c) b
\eeq 
be the translate of the set $\varphi_J(A)$ by the subset sum $s(J^c) b$.  
We have $\Phi_{\emptyset}(A,b) =  \left( \sum_{j=1}^h  c_j \right) b$ 
and $\Phi_J(A,b) = \varphi(A)$ if $J = \{1,\ldots, h\}$.

\bl                                                                   \label{LinearPerturb:lemma:disjointSidon}
Let $\varphi = \sum_{i=1}^h c_i x_i$ be a linear form with coefficients in the field $\mbF$.  
Let $V$ be a vector space over \mbF.  For every subset $A$ of $V$ and $b \in V \setminus A$, 
\beq                                                                              \label{LinearPerturb:phi-partition}
\varphi\left( A \cup \{b\} \right) = \bigcup_{J \subseteq \{1,\ldots, h\} } \Phi_J(A,b).  
\eeq
If  $A \cup \{b\}$ is a $\varphi$-Sidon set, then 
\beq                                                                              \label{LinearPerturb:phi-J}
\left\{ \Phi_J(A,b): J \subseteq  \{1,\ldots, h\} \right\}
\eeq
is a set of pairwise disjoint sets.  

If $A$ is a $\varphi$-Sidon set and~\eqref{LinearPerturb:phi-J}  
is a set of pairwise disjoint sets, then $A \cup \{b\}$ is a $\varphi$-Sidon set.  
\el

\begin{proof}
 If $w \in \varphi\left( A \cup \{b\} \right)$, then there exist 
 vectors $v_1,\ldots, v_h \in A \cup \{b\}$ 
 such that 
 \[
w = \varphi(v_1,\ldots, v_h)  = \sum_{i=1}^h c_i v_i.
\]
Let $J = \{j \in \{1,\ldots, h\}: v_j = a_j \in A$.  We have $J^c = \{j \in \{1,\ldots, h\}: v_j = b\}$ 
and 
\begin{align*}
w & = \sum_{i=1}^h c_i v_i = \sum_{j \in J} c_j a_j + \sum_{j\in J^c} c_j b 
\in \varphi_J(A) + s(J^c) b = \Phi_J(A,b). 
\end{align*} 
Conversely, if $w \in  \Phi_J(A,b)$ for some $J = \{j \in \{1,\ldots, h\}$, 
then there exist $a_j \in A$ for all $j \in J$ such that 
\[
w  = \sum_{j \in J} c_j a_j + \sum_{j\in J^c} c_j b \in \varphi\left( A \cup \{b\} \right).
\]
This proves~\eqref{LinearPerturb:phi-partition}.   
It follows that if $A \cup \{b\}$ is a $\varphi$-Sidon set, then~\eqref{LinearPerturb:phi-J}  
is a set of pairwise disjoint sets.  

Suppose that $A$ is a Sidon set and that the sets $\Phi_{J}(A,b)$ are pairwise disjoint  
for all  $J \subseteq \{1,\ldots, h\}$. 
Let $u_1,\ldots, u_h, v_1,\ldots, v_h  \in A \cup \{b\}$.    
Consider the sets  
\[
J_1 = \{ j \in \{1,\ldots, h\} : u_j \neq  b\} \qqand J_2 = \{ j \in \{1,\ldots, h\} : v_j \neq  b\}  
\]
and the complementary sets  
\[
J_1^c = \{ j \in \{1,\ldots, h\} : u_j =  b\} \qqand J_2^c = \{ j \in \{1,\ldots, h\} : v_j =  b\}. 
\]
We have 
\[
\varphi(u_1,\ldots, u_h) = \sum_{j \in J_1} c_ju_j +  \left( \sum_{j \in J_1^c} c_j \right) b
\in \Phi_{J_1}(A,b) 
\]
and
\[
\varphi(v_1,\ldots, v_h) = \sum_{j \in J_2} c_jv_j + \left( \sum_{j \in J_2^c} c_j \right) b  
\in \Phi_{J_2}(A,b). 
\]

If $J_1 \neq J_2$, then $\Phi_{J_1}(A,b)  \cap \Phi_{J_2}(A,b)  = \emptyset$ and 
$\varphi(u_1,\ldots, u_h) \neq \varphi(v_1,\ldots, v_h)$.  

If $J_1 = J_2 = \emptyset$, then $(u_1,\ldots, u_h) = (b,\ldots, b) = (v_1,\ldots, v_h)$.  

If $J_1 = J_2 \neq \emptyset$, then $J_1^c = J_2^c$ and 
\[
\sum_{j \in J_1^c} c_j =  \sum_{j \in J_2^c} c_j.   
\]
It follows that  
\[
 \sum_{j \in J_1} c_ju_j  = \sum_{j \in J_1} c_jv_j.
\]
Because $A$ is a $\varphi_{J_1}$-Sidon set, we have 
$u_j = v_j$ for all $j \in J_1$, hence 
$u_i = v_i$ for all $i \in \{1,\ldots, h\}$.  
Thus, if $A$ is a Sidon set and the sets $\Phi_{J}(A,b)$ are pairwise disjoint, 
then $A \cup \{ b\}$ is a $\varphi$-Sidon set.  
This completes the proof.  
\end{proof}

\bl                                                                               \label{LinearPerturb:lemma:disjointPhi-0}
Let $\varphi = \sum_{i=1}^h c_i x_i$ be a linear form with coefficients 
in the field $\mbF$.  
Let $V$ be a vector space over \mbF, let $X$ be an infinite subset of $V$, 
and let $B$ be a finite subset of $X$.  
If the  linear form $\varphi$ has property $N$,  
then there exists $b \in X$ such that, for all subsets $J$ of $\{1,\ldots, h\}$, 
the sets $\Phi_J(B,b)$ are pairwise disjoint.  
\el

\begin{proof}
Let $J_1$ and $J_2$ be distinct subsets of $\{1,\ldots, h\}$.  
For all $x \in X$, we have 
\beq                                                      \label{LinearPerturb:phiIntersect-0}
\Phi_{J_1}(B,x) \cap \Phi_{J_2}(B,x) \neq \emptyset
\eeq
if and only if there exist elements $b_{1,j} \in B$ for all $j \in J_1$  
and  $b_{2,j} \in B$ for all $j \in J_2$ such that 
\beq                                                      \label{LinearPerturb:phiIntersect-3}
 \sum_{j \in J_1} c_jb_{1,j} + \left( \sum_{j \in J^c_1} c_j \right) x
 =  \sum_{j \in J_2} c_j b_{2,j} + \left( \sum_{j \in J^c_2} c_j\right) x. 
\eeq
Let $K = J^c_1 \cap J^c_2$.     
The sets $ I_1 = J^c_1 \setminus K$ and $I_2 = J^c_2 \setminus K$ are disjoint. 
If $I_1 = I_2 = \emptyset$, then $J^c_1 = K = J^c_2$ and  
$J_1  = J_2$, which is absurd.  
Therefore, $I_1$ and $I_2$ are disjoint sets, not both empty.

Because the linear form $\varphi$ has property $N$, we have 
\[
\sum_{j \in I_1} c_j  \neq \sum_{j \in I_2} c_j 
\]
and so  
\[
c = \sum_{j \in I_2} c_j  - \sum_{j \in I_1} c_j \neq 0.
\]
Thus, $c \in \mbF\setminus \{ 0\}$ and so $c$  is invertible in $\mbF$.  
From~\eqref{LinearPerturb:phiIntersect-3} we obtain 
\begin{align*}
 \sum_{j \in J_1} c_j b_{1,j} -  \sum_{j \in J_2} c_j b_{2,j} 
& =  \left( \sum_{j \in J^c_2} c_j \right)x - \left( \sum_{j \in J^c_1} c_j \right)x \\
& =  \left(  \sum_{j \in I_2} c_j - \sum_{j \in I_1} c_j \right) x \\
& =  c x 
\end{align*} 
and so 
\beq                                                      \label{LinearPerturb:xx}
x = c^{-1} 
\left(  \sum_{j \in I_1} c_j b_{1,j} -  \sum_{j \in I_2} c_j b_{2,j}  \right).  
\eeq
Because the set $B$ is finite, the set $B'$ of elements in $X$ 
of the form~\eqref{LinearPerturb:xx} is also finite.
Because the set $X$ is infinite,  the set $X \setminus (B \cup B')$ is infinite.  
For all $b \in X \setminus (B \cup B')$, the set  
 $\{\Phi_J(B,b) : J \subseteq \{1,\ldots, h\} \}$ consists of pairwise disjoint sets.  
This completes the proof.  
\end{proof}

\bt                                                                \label{LinearPerturb:theorem:exist-V}
Let \mbF\ be a field, let $V$ be an infinite vector space over the field \mbF, 
and let $X$ be an infinite subset of $V$.   
Let $\varphi(x_1,\ldots, x_h) = \sum_{i=1}^h c_ix_i$ be a linear form 
with nonzero coefficients $c_i \in \mbF$.   
The following are equivalent:
\benum
\item[(i)] 
The set $X$ contains an infinite $\varphi$-Sidon set $A$.
\item[(ii)]
The set $X$ contains a $\varphi$-Sidon set $A$ with $|A| \geq 2$. 
\item[(iii)]
The linear form $\varphi$ has property $N$.  
\eenum
\et

\begin{proof}
Condition~(i) implies~(ii).  
It was proved in Section~\ref{LinearPerturb:section:forms} that~(ii) implies~(iii).  
We shall prove that (iii) implies (i).

Suppose that the  linear form $\varphi$ has property $N$. 
We construct  inductively an infinite $\varphi$-Sidon set $A$ contained in $X$.
For all $a_1 \in X$, the set $A_1 = \{a_1\}$ is a $\varphi$-Sidon set, 
because every set with one element is $\varphi$-Sidon.    
Let $A_n = \{a_1,\ldots, a_n\}$ be a $\varphi$-Sidon set $A$ contained in $X$.  
By Lemma~\ref{LinearPerturb:lemma:disjointPhi-0}, there exists $a_{n+1} \in X$ 
such that 
\[
\Phi_{J_1}(A_n,a_{n+1}) \cap \Phi_{J_2}(A_n,a_{n+1}) = \emptyset
\]
if  $J_1$ and $J_2$ are distinct subsets of $\{1,\ldots, h\}$.  
It follows from Lemma~\ref{LinearPerturb:lemma:disjointSidon} 
that the set $A_{n+1} = A_n \cup \{a_{n+1} \}$ is a $\varphi$-Sidon set.  
This completes the proof.  
\end{proof}

\section{Perturbations of linear forms}                    \label{LinearPerturb:section:Perturbations} 

An absolute value  on a field \mbF\ is a function $| \ |: \mbF  \rightarrow \R$ such that 
\benum
\item[(i)]
$|c| \geq 0$ for all $c\in \mbF$, and $|c|=0$ if and only if $c=0$,
\item[(ii)]
$|c_1c_2| = |c_1| \ |c_2|$ for all $c_1,c_2 \in \mbF$,  
\item[(ii)]
$|c_1+c_2| \leq |c_1| + |c_2|$ for all $c_1,c_2 \in \mbF$.
\eenum
The absolute value $| \ |$ on \mbF\ is \emph{trivial} if $|c|=1$ for all $c \neq 0$, 
and  \emph{nontrivial} if $|c| \neq 1$ for some $c \neq 0$. 
The usual absolute value and the $p$-adic absolute values are the nontrivial absolute values on \Q.

Let $V$ be a vector space over \mbF.  A \emph{norm} on $V$ with respect to an absolute value 
$| \ |$ on \mbF\ is a function $\| \ \|:V\rightarrow \R$ such that 
\benum
\item[(i)]
$\| v \| \geq 0$ for all $v\in V$, and $\| v\| = 0$ if and only if $v=0$,
\item[(ii)]
$\|cv\| = |c| \ \|v\|$ for all $c\in \mbF$ and $v\in V$,
\item[(iii)]
$\|v+w\| \leq \| v \| + \| w\|$ for all $v,w \in V$.
\eenum 

For example, if $| \ |$ is an absolute value on \mbF\ and $V = \mbF^n$, then, 
for every vector $\mbx = \vectorsmallxn \in V$, the functions 
\[
\| \mbx \|_1 = \sum_{j=1}^n |x_j|
\]
and
\[
\|\mbx\|_{\infty} = \max\{ |x_j|: j = 1,\ldots, n\}
\]
are norms on $V$ with respect to $| \ |$.

If $| \ |$ is a nontrivial absolute value on \mbF, then there exists $c \in \mbF$ with $|c| \neq 0$
and $|c| \neq 1$.  
If $|c| > 1$, then $0 < |1/c| = 1/|c| < 1$.  
If $0 < |c_0|<1$, then  
\[
0 < \left| c_0^{n+1}\right| = \left| c_0\right|^{n+1} < \left| c_0\right| ^n = \left| c_0^n \right| 
\]
 for all $n \in \N$.  
Thus, the field \mbF\ is infinite and 
\beq                                                      \label{LinearPerturb:infiniteF}
\inf\{ |c|: c \in \mbF\setminus \{0\} \} =\inf\{ |c_0^n|: n=1,2,3,\ldots\} = 0.
\eeq

Let $V$ be a nonzero normed vector space with respect to a nontrivial absolute value on the field \mbF.  
Let $v_0 \in V\setminus \{ 0\}$.  Let $c_0 \in \mbF$ with $0 < |c_0|<1$.  
For all $n \in \N$ we have $c_0^n v_0 \neq 0$ and  
\[
0 < \left\| c_0^{n+1} v_0 \right\| =   \left| c_0^{n+1} \right|   \left\| v_0 \right\| 
<  \left| c_0^{n} \right|  \left\| v_0 \right\| = \left\| c_0^n v_0 \right\| 
\]
Thus, the vector space $V$ is infinite and 
\beq                                                      \label{LinearPerturb:infiniteV}
\inf\{ |x|: x \in V \setminus \{0\} \} =\inf\{ |c_0^n v_0|: n=1,2,3,\ldots\} = 0.
\eeq

\bl                                                                               \label{LinearPerturb:lemma:disjointPhi}
Let \mbF\ be a field with a nontrivial absolute value.  
Let $V$ be a nonzero vector space over \mbF\
that has a norm with respect to the absolute value on \mbF.  
Let $A'$ be a  finite subset of $V$ and let $b \in V$.  

Let $\varphi= \sum_{i=1}^h c_i x_i$ 
be a linear form with coefficients $c_i \in \mbF$.   
If the  linear form $\varphi$ has property $N$,  
then for every $\varepsilon > 0$ there are infinitely many nonzero vectors $a \in V$ such that 
\[
\|a-b \| < \varepsilon
\]
and, for all subsets $J$ of $\{1,\ldots, h\}$, 
the sets 
\[     
\Phi_J(A',a) = \varphi_J(A') + \left( \sum_{j \in J^c} c_j\right) a 
\]
are pairwise disjoint.  
\el

\begin{proof}
If $A' = \emptyset$, then $\varphi_J(A') = \{ 0 \}$ 
for all $J \subseteq \{1,\ldots, h\}$ 
and $\Phi_J(A',a) = \left\{  \left( \sum_{j \in J^c} c_j\right) a \right\}$.
Because $\varphi$ has property $N$,  for every nonzero vector $a \in V$ 
the vectors $ \left( \sum_{j \in J^c} c_j\right) a  = s(J^c) a$ are distinct and so the sets 
$\Phi_J(A',a)$ are pairwise disjoint.  
Choose any of the infinitely many nonzero vectors $a$ such that $\|a-b\|<\varepsilon$.   

Let $A' \neq \emptyset$ and $x \in V$.  
For distinct subsets $J_1$ and $J_2$ of $\{1,\ldots, h\}$, we have 
\beq                                                      \label{LinearPerturb:phiIntersect}
\Phi_{J_1}(A',b+x) \cap \Phi_{J_2}(A',b+x) \neq \emptyset
\eeq
if and only if there exist vectors $a_{1,j} \in A'$ for all $j \in J_1$  
and  $a_{2,j} \in A'$ for all $j \in J_2$ such that 
\beq                                                      \label{LinearPerturb:phiIntersect-2}
 \sum_{j \in J_1} c_ja_{1,j} + \sum_{j \in J^c_1} c_j(b+x) 
 =  \sum_{j \in J_2} c_j a_{2,j} + \sum_{j \in J^c_2} c_j(b+x). 
\eeq
Let $K = J^c_1 \cap J^c_2$.  The sets    
$ I_1 = J^c_1 \setminus K$ and $I_2 = J^c_2 \setminus K$ are disjoint.   
If $I_1 = I_2 = \emptyset$, then $K = J^c_1  = J^c_2$ and  so 
$J_1  = J_2$, which is absurd.  Therefore, 
the sets $I_1$ and $I_2$ are disjoint sets, not both empty. 

Because the linear form $\varphi$ has property $N$, we have 
\[
\sum_{j \in I_1} c_j  = s(I_1) \neq s(I_2) =  \sum_{j \in I_2} c_j 
\]
and so  
\[
c = \sum_{j \in I_2} c_j  - \sum_{j \in I_1} c_j \neq 0.
\]
Thus, the scalar $c$  is invertible in $\mbF$.  
From~\eqref{LinearPerturb:phiIntersect-2} we obtain 
\begin{align*}
 \sum_{j \in J_1} c_j a_{1,j} -  \sum_{j \in J_2} c_j a_{2,j} 
& =  \sum_{j \in J^c_2} c_j(b+x) - \sum_{j \in J^c_1} c_j(b+x) \\
& =  \sum_{j \in I_2} c_j(b+x) - \sum_{j \in I_1} c_j(b+x) \\
& =  \left(  \sum_{j \in I_2} c_j - \sum_{j \in I_1} c_j \right) (b+x) \\
& =  c (b+x) 
\end{align*} 
and 
\beq                                                      \label{LinearPerturb:x}
x = c^{-1} 
\left(  \sum_{j \in J_1} c_j a_{1,j} -  \sum_{j \in J_2} c_j a_{2,j}  \right) - b.
\eeq
Because the set $A'$ is nonempty and  finite, the set $X$ of vectors $x$  in $V$ 
of the form~\eqref{LinearPerturb:x} is also nonempty and finite.  
If $X = \{0\}$, let $\delta = 1$.  
If $X \neq \{0\}$, let 
\beq                                                      \label{LinearPerturb:x-ineq-1}
\delta = \min\{ \| x \|:x\in X\setminus \{0\} \} > 0
\eeq
and let 
\beq                                                      \label{LinearPerturb:x-ineq-2}
\varepsilon' =  \min(\delta,\varepsilon) > 0. 
\eeq
By~\eqref{LinearPerturb:infiniteV}, 
there are infinitely many vectors $x_0$ in $V$ such that 
\beq                                                      \label{LinearPerturb:x-ineq}
0 < \| x_0 \| < \varepsilon'.
\eeq
It follows from~\eqref{LinearPerturb:x-ineq-1} and~\eqref{LinearPerturb:x-ineq-2}  
that each such vector satisfies $x_0 \notin X$,  and so 
\[
\Phi_{J_1}(A',b+x_0) \cap \Phi_{J_2}(A',b+x_0) = \emptyset
\]
for all distinct subsets $J_1$ and $J_2$ of $\{1,\ldots, h\}$.  
Choosing $a = b +x_0$ completes the proof.  
\end{proof}

Let \mbF\ be a field with a nontrivial absolute value, and let $V$ be a 
vector space over \mbF\ that has a norm with respect to the absolute value on \mbF.  
Let $\N= \{1,2,3,\ldots \}$ be the set of positive integers.  
Let $A = \{a_k : k \in \N\}$ and $B = \{b_k : k \in \N \}$ 
be sets of not necessarily distinct vectors in $V$. 
Let $\varepsilon = \{\varepsilon_k : k \in \N \}$ be a set of positive real numbers.  
The set  $B$  is an  
\emph{$\varepsilon$-perturbation}\index{perturbation}  
of the set $A$ if 
\[
\|a_k - b_k \| < \varepsilon_k
\]
for all $k \in \N$.

\bt                                                                       \label{LinearPerturb:theorem:asymptotic}
Let \mbF\ be a field  with a nontrivial absolute value  
and let $V$ be a vector space over \mbF\ that has a norm with respect to the absolute value 
on \mbF.  Let $\varepsilon = \{\varepsilon_k: k = 1,2,3,\ldots \}$ be a set
of positive real numbers.
Let $\varphi$ be a linear form with coefficients in \mbF\ that has property $N$.  
For every set $B = \{ b_k: k = 1,2,3,\ldots \}$ of vectors in $V$, 
there is a $\varphi$-Sidon set $A= \{a_k: k = 1,2,3,\ldots \}$ of vectors in $V$ 
such that 
 \beq                                                                                   \label{LinearPerturb:inequality}
\| a_k - b_k \| < \varepsilon_k 
\eeq
 for all $k = 1,2,3,\ldots$. 
\et

\begin{proof}
We construct the set $A$  inductively.  
Begin by choosing $a_1 = b_1$.  
Every set with one element is a $\varphi$-Sidon set, 
and so $A_1 = \{a_1\}$ is a $\varphi$-Sidon set 
such that $\| a_1-b_1 \| = 0 < \varepsilon_1$.  

Let $n \geq 1$, and let $A_n = \{a_1,\ldots, a_n \}$ be a $\varphi$-Sidon set 
that satisfies inequality~\eqref{LinearPerturb:inequality}  for all  $k \in \{1,\ldots,n\}$.  
Applying Lemma~\ref{LinearPerturb:lemma:disjointPhi} to the finite set $A' = A_n$ 
and the vector $b = b_{n+1}$, we obtain a vector $a_{n+1} \in V$ such that 
$\|a_{n+1} - b_{n+1}\| < \varepsilon_{n+1}$  and the sets 
$\Phi_J(A_n, a_{n+1})$ are pairwise disjoint for all $J \subseteq \{1,\ldots, h\}$.  
The set $A_n$ is $\varphi$-Sidon, and so, by Lemma~\ref{LinearPerturb:lemma:disjointSidon}, 
the set $A_{n+1} = A_n \cup \{ a_{n+1} \}$ is a  $\varphi$-Sidon set.
This completes the proof.  
\end{proof}

\bt                                                                \label{LinearPerturb:theorem:limit}
Let \mbF\ be a field  with a nontrivial absolute value, 
and let $\varphi$ be a linear form with coefficients in \mbF\ that has property $N$.  
Let $V$ be a vector space over \mbF\ that has a norm with respect to absolute value on \mbF. 
For every set $B = \{b_k: k=1,2,3,\ldots \}$ of vectors in $V$, there exists a $\varphi$-Sidon set
 $A = \{a_k: k=1,2,3,\ldots \}$ in $V$ such that 
 \[
\lim_{k \rightarrow \infty} \| a_k  - b_k \| = 0.
 \]
\et

\begin{proof}
This follows fromTheorem~\ref{LinearPerturb:theorem:asymptotic} 
applied to any sequence $\varepsilon = \{\varepsilon_k: k = 1,2,3,\ldots \}$ of positive numbers such that 
$\lim_{k\rightarrow \infty} \varepsilon_k = 0$.  
\end{proof}

\section{$p$-adic $\varphi$-Sidon sets} 

Let $\PP = \{2,3,5,\ldots \}$ be the set of prime numbers.  
For every prime number $p$, let $| \ |_p$ be the usual $p$-adic absolute value with $|p|_p = 1/p$.
Every integer $r$ satisfies $|r|_p \leq 1$. 

\bl                                                                               \label{LinearPerturb:lemma:p-adic}
Let $\varphi= \sum_{i=1}^h c_i x_i$ 
be a linear form with rational coefficients $c_i$ that satisfies property $N$.   
Let $\PP_0$ be a nonempty finite set of prime numbers.  
Let $A'$ be a  finite set of integers and let $b$ be an integer.  
For every $\varepsilon > 0$ there are infinitely many positive integers $a$ such that 
\[
|a-b|_{p} < \varepsilon
\]
for all $p \in \PP_0$ and the sets 
\[     
\Phi_J(A',a) = \varphi_J(A') + \left( \sum_{j \in J^c} c_j\right) a 
\]
are pairwise disjoint for all subsets $J$ of $\{1,\ldots, h\}$.    
\el

\begin{proof}
Let $\varepsilon' > 0$.  Choose a positive integer $k$ such that 
\[
\frac{1}{2^k} < \varepsilon'.
\]
The integer $b$ is not necessarily positive, 
but for all sufficiently large positive integers $r$ we have 
\beq                                                                           \label{LinearPerturb:p-adic}
a=b+r \prod_{p \in \PP_0} p^k > 0.
\eeq 
Let $p \in \PP_0$.  
For all integers $r$ satisfying~\eqref{LinearPerturb:p-adic} we have   
\[
|a-b|_{p}  = |r|_{p} \prod_{p \in \PP_0} |p^k|_p \leq   |p^k|_{p}   
= \frac{1}{p^k}  \leq \frac{1}{2^k} < \varepsilon' .
\]

The proof of Lemma~\ref{LinearPerturb:lemma:p-adic} is the same as the proof of 
Lemma~\ref{LinearPerturb:lemma:disjointPhi} until the choice of $x_0$,   
at which point we choose a positive integer $x_0 = r \prod_{p \in \PP_0} p^k$ 
that satisfies inequality~\eqref{LinearPerturb:p-adic}.  
This completes the proof.  
\end{proof}

\bt                                     \label{LinearPerturb:theorem:Positive}
Let $\varphi$ be a linear form with rational coefficients that satisfies property $N$.  
Let $\{\varepsilon_k : k=1,2,3,\ldots \}$ be a sequence of positive real numbers 
and let $\{p_k:k=1,2,3,\ldots\}$ be a sequence of prime numbers.    
For every sequence of integers $B = \{b_k:k=1,2,3,\ldots\}$, 
there exists a strictly increasing sequence of positive integers $A =\{ a_k:k=1,2,3,\ldots\}$ such that 
$A$ is a $\varphi$-Sidon set and 
\[
 | a_k - b_k |_{p_j} < \varepsilon_k 
\]
for all $k \in \N$ and $j \in \{1,\ldots, k\}$.  
\et

\begin{proof}
The proof  of Theorem~\ref{LinearPerturb:theorem:Positive} is an inductive construction 
based on Lemma~\ref{LinearPerturb:lemma:p-adic}.
Choose a positive integer $k_1$ such that $1/p_1^{k_1} < \varepsilon_1$ 
and $b_1+p_1^{k_1} > 0$.
Let $a_1 = b_1+p_1^{k_1}$.  The set $A_1 = \{a_1\}$ is a $\varphi$-Sidon set and 
 $ | a_1 - b_1 |_{p_1} < \varepsilon_1$.

For $n \geq 1$, let $A_n = \{a_1,\ldots, a_n\}$ be a set of  positive integers 
with $a_1 < \cdots < a_n$ 
such that $A_n$ is a $\varphi$-Sidon set and 
\[
 | a_k - b_k |_{p_j} < \varepsilon_k 
\]
for all $k \in \{1,\ldots, n\}$ and $j \in \{1,\ldots, k \}$.  
We apply Lemma~\ref{LinearPerturb:lemma:p-adic} 
to the set $A' =A_n$, the integer $b = b_{n+1}$,  the finite set of primes 
$\PP_0=\{p_1,\ldots, p_n,p_{n+1}\}$, and 
$\varepsilon' = \varepsilon_{n+1} > 0$ 
to obtain an integer $a_{n+1} > a_n$ such that 
\[
|a_{n+1} - b_{n+1} |_{p_j} <  \varepsilon_{n+1}
\] 
for all $j \in \{1,\ldots, n,n+1\}$ and the sets $\Phi_J(A_{n+1},a _{n+1})$ 
are pairwise disjoint for all $J \subseteq \{1,\ldots, h\}$.  
It follows from Lemma~\ref{LinearPerturb:lemma:disjointSidon} that 
$A_{n+1}$ is a $\varphi$-Sidon set.  
This completes the proof.  
\end{proof}

\bt                                    \label{LinearPerturb:theorem:Positive-2}
Let $\varphi$ be a linear form with  rational coefficients that satisfies property $N$.   
Let $B = \{b_k : k=1,2,3,\ldots \}$ be a sequence of integers.  
There exists a strictly increasing $\varphi$-Sidon set of positive integers $A =\{ a_k : k=1,2,3,\ldots \}$ 
such that, for every prime number $p$, the set $A$ is $p$-adically 
asymptotic to $B$ in the sense that 
\[
\lim_{k \rightarrow \infty} | a_k - b_k |_p = 0. 
\]
\et

\begin{proof}
This follows from Theorem~\ref{LinearPerturb:theorem:Positive} 
applied to the set of all prime numbers and 
any sequence $\varepsilon = \{\varepsilon_k: k = 1, 2, 3, \ldots \}$ of positive numbers such that 
$\lim_{k\rightarrow \infty} \varepsilon_k = 0$.  
\end{proof}

\section{Growth of $\varphi$-Sidon sets} 

Let $f(t)$  be a real-valued  or complex-valued function defined for $t \geq t_0$.
Let $g(t)$ be  positive function defined for $t \geq t_0$.
We write 
\[
f(t) \ll g(t)  
\]
if there  exist constants $C_1 > 0$ and $t_1\geq t_0$ such that $|f(t)| \leq C_1 g(t)$ 
for all $t \geq t_1$.  
We write 
\[
f(t) \gg g(t)  
\]
if there exist constants $C_2 > 0$ and $t_2 \geq t_0$ such that $|f(t)| \geq C_2 g(t)$ 
for all $t \geq t_2$.

Let  $A$ be a set of positive integers.   
The \emph{growth function} or \emph{counting function} of $A$ is the function $A(n)$ 
that counts the number of positive integers in the set $A \cap \{1,\ldots, n\} $.  
The number of $h$-fold sums of integers taken from the set $A \cap \{1,\ldots, n\} $ is 
\[
\binom{A(n)+h-1}{h}
\]
and each of these sums is at most $hn$.  
If $A$ is a classical $h$-Sidon set, then these sums are distinct and 
\[
\frac{A(n)^h}{h!} \leq \binom{A(n)+h-1}{h} \leq hn
\]
This  simple counting argument proves that 
\[
A(n) \ll n^{1/h}.  
\]
The upper bound is tight.  
Bose and Chowla~\cite{bose-chow62} proved that 
for every positive integer $n$ there exist finite  Sidon sets $A$ with 
\[
A \subseteq \{1,\ldots, n\} \qqand \card(A) \gg n^{1/h}.
\]

We do not have best possible upper bounds for infinite Sidon sets.  
Erd\H os (in St{\"o}hr~\cite{stoh55}) constructed an infinite Sidon set $A$ of order 2 with 
\[
\limsup_{n\rightarrow\infty} \frac {A(n)}{\sqrt{n}} \geq \frac{1}{2}
\]
and so $A(n) \gg \sqrt{n}$ for infinitely many $n$, 
but he also proved that every classical Sidon set of order 2 satisfies 
\[
\liminf_{n\rightarrow\infty} A(n) \sqrt{ \frac {\log n}{n}} \ll 1
\]
and so $A(n) \ll \sqrt{n/\log n}$ 
for infinitely many $n$.

It is of interest to obtain upper bounds for the size of $\varphi$-Sidon sets. 
Let $\mbF$ be a field with an absolute value.  
The \emph{counting function} of  a subset $X$ of \mbF\ is
\[
X(t) = \card \left( x\in X: |x| \leq t \right).  
\] 

\bt                     \label{LinearPerturb:theorem:growth}
Let $\mbF$ be a field with an absolute value.  
Let $\varphi = \sum_{i=1}^h c_ix_i$ be a  linear form with coefficients in $\mbF$,
and let $C = \sum_{i=1}^h |c_i|$.  
Let $X$ be a subset of \mbF\ such that $\varphi(X) \subseteq X$.  
If $A$ is a $\varphi$-Sidon subset of $X$, then   
\[
A(t) \leq X(Ct)^{1/h} 
\]
for all $t \geq 0$.  
\et

\begin{proof}
Let $A' = \{a \in A: |a| \leq t\}$.  
We have $A(t) = \card(A')$ and, because $A$ is a  $\varphi$-Sidon set,  
\[
A(t)^h = \card(\varphi(A')).
\] 
If $a_1,\ldots, a_h \in A'$, then $b = \varphi(a_1,\ldots, a_h)  \in \varphi(A') \subseteq X$ and  
\begin{align*}
|b| & = \left| \varphi(a_1,\ldots, a_h) \right|  = \left| \sum_{i=1}^h c_ia_i \right| \\
& \leq \sum_{i=1}^h \left|  c_ia_i \right|  \leq  \sum_{i=1}^h \left|  c_i  \right|  \max(|a_i|:i=1,\ldots h) \\ 
& \leq C t. 
\end{align*}
Therefore, 
\[
A(t)^h = \card(\varphi(A')) \leq \card\{ x\in X: |x| \leq Ct\} = X(Ct)  
\]
and   
\[
A(t) \leq X(Ct)^{1/h}. 
\] 
This completes the proof.  
\end{proof}

Let $\varphi = \sum_{i=1}^h c_i x_i$ be a  linear form with nonzero rational coefficients.  
Let $m$ be a common multiple of the the denominators of the coefficients $c_1,\ldots, c_h$, 
and let $d$ be the greatest common divisor of the integers $mc_1, \ldots, mc_h$.   
Let $c'_i = mc_i/d$  for $i \in \{1,\ldots, h\}$. 
The integers $c'_i = mc_i/d$ are nonzero and relatively prime.   
Consider the linear form $\varphi' = \sum_{i=1}^h c'_i x_i$.  We have 
\[
\varphi 
= \frac{d}{m} \sum_{i=1}^h \frac{mc_i}{d} x_i = \frac{d}{m} \sum_{i=1}^h c'_i x_i 
= \frac{d}{m} \varphi'.   
\]
It follows that a set is a $\varphi$-Sidon set if and only if it is a $\varphi'$-Sidon set.  
Thus, in the study of $\varphi$-Sidon sets, a linear form with nonzero rational coefficients 
can be replaced with a linear form with nonzero relatively prime integer coefficients.

\bt
Let $\varphi = \sum_{i=1}^h c_ix_i$ be a linear form with integer coefficients.   
If $A$ is a $\varphi$-Sidon set of integers, then   
\[
A(t) = \{a\in A: |a| \leq t \} \ll t^{1/h}.   
\] 
\et

\begin{proof}
We have $\varphi(\Z) \subseteq \Z$.  
Let $[t]$ denote the integer part of the real number $t$.
With the usual absolute value, the counting function of \Z\ is $\Z(t) = 2[t]+1 \leq 2t+1$.  
Applying Theorem~\ref{LinearPerturb:theorem:growth} with $X=\Z$, we obtain 
\[
A(t) \leq \Z(Ct)^{1/h} \leq (2Ct+1)^{1/h} \ll t^{1/h}.
\]
This completes the proof.  
\end{proof}

 \bt                \label{LinearPerturb:theorem:LowerGrowth}
Let $\varphi = \sum_{i=1}^h c_ix_i$ be a linear form with integer coefficients 
that satisfies condition $N$.    
There exists an infinite $\varphi$-Sidon set $A = \{a_k: k \in \N\}$ 
of distinct  positive integers such that 
\beq                                         \label{LinearPerturb:GrowthIneq}
a_{k+1}   \leq 4^h k^{2h-1} + k 
\eeq
for all $k \in \N$.   
\et

\begin{proof}
We construct the $\varphi$-Sidon set $A = \{a_k: k \in \N\}$ inductively.  
The set $A_1 = \{a_1\}$ is a $\varphi$-Sidon set for every integer $a_1$.  
Let $a_1 = 1$.  

Let $k \geq 1$ and let $A_{k} = \{a_1,\ldots, a_{k} \}$ be a  $\varphi$-Sidon set of positive integers.
Let $b$ be a positive integer.  
By Lemma~\ref{LinearPerturb:lemma:disjointSidon}, the set $A_{k} \cup \{b\}$ 
is a $\varphi$-Sidon set if and only if the sets 
\[
\Phi_J(A_{k},b) = \varphi_J(A_{k}) +  \left( \sum_{j \in J^c} c_j\right) b
\] 
are pairwise disjoint for all $J \subseteq \{1,\ldots, h\}$.  

Let $J_1$ and $J_2$ be distinct subsets of $\{1,\ldots, h\}$.  
The sets $J_1\setminus (J_1\cap J_2)$ and $J_2 \setminus (J_1\cap J_2)$ 
are distinct and disjoint.  
We have 
\[
\Phi_{J_1}(A_{k},b) \cap \Phi_{J_2}(A_{k},b) \neq \emptyset 
\]
if and only if there exist integers $a_{1,j} \in A_{k}$ for all $j \in J_1$ 
and $a_{2,j} \in A_{k}$ for all $j \in J_2$ such that 
\beq                                                 \label{LinearPerturb:bc}
 \sum_{j\in J_1} c_ja_{1,j} +  \left( \sum_{j \in J^c_1} c_j\right)b 
 = \sum_{j\in J_2} c_ja_{2,j} +  \left( \sum_{j \in J^c_2} c_j\right) b. 
\eeq
The integer 
\begin{align*} 
c & =   \sum_{j \in J^c_2} c_j - \sum_{j \in J^c_1} c_j = s(J_2^c) - s(J_1^c) \\\
& = s\left( J_1\setminus (J_1\cap J_2) \right) - s \left(J_2 \setminus (J_1\cap J_2) \right) 
\end{align*}
 is nonzero because the linear form $\varphi$ satisfies condition $N$.  
The integer $b$ satisfies equation~\eqref{LinearPerturb:bc} 
if and only if 
 \beq                                         \label{LinearPerturb:GrowthEquation}
cb = \sum_{j\in J_1} c_j a_{1,j} - \sum_{j\in J_2} c_j a_{2,j}.    
\eeq
Thus, there is at most one integer $b$ that satisfies equation~\eqref{LinearPerturb:GrowthEquation}.   

Let $\card(J_1) = j_1$ and $\card(J_2) = j_2$.  
The sets $J_1$ and $J_2$ are distinct subsets of $\{1,\ldots, h\}$ 
and so 
\[
j_1 + j_2 \leq 2h-1.  
\] 
The number of integers of the form 
\[
 \sum_{j\in J_1 } c_j a_{1,j}- \sum_{j\in J_2} c_j a_{2,j}
 \]
 with $a_{1,j} \in A_{k}$ and $a_{2,j} \in A_{k}$ is at most ${k}^{j_1+j_2}$.  
 The number of ordered pairs $(J_1,J_2)$ of  subsets of $\{1,\ldots, h\}$ 
 of cardinalities $j_1$ and $j_2$, respectively,  is 
 \[
 \binom{h}{j_1}  \binom{h}{j_2}.
  \]
 Thus, the number of equations of the form~\eqref{LinearPerturb:GrowthEquation} 
is at most 
\begin{align*}
\underbrace{\sum_{j_1=0}^h  \sum_{j_2=0}^{h} }_{j_1+j_2\leq 2h-1}  \binom{h}{j_1} \binom{h}{j_2}  {k}^{j_1+j_2} 
&  \leq \sum_{j_1=0}^h  \binom{h}{j_1}  \sum_{j_2=0}^{h} \binom{h}{j_2} {k}^{2h-1} \\
&= 4^h k^{2h-1} 
\end{align*}
and so there are at most $4^h k^{2h-1} + k$ positive integers $b$ such that $b \notin A_{k}$ 
and $A_{k} \cup \{b\}$ is not a $\varphi$-Sidon set.   
It follows that there exists a positive integer $a_{k +1}$ such that
\benum
\item [(i)] 
$a_{k +1} \notin A_{k}$, 
\item [(ii)] 
$A_{k +1}= A_{k} \cup \{ a_{k +1} \}$ is a $\varphi$-Sidon set, 
\item [(iii)] 
 $a_{k +1} \leq 4^h k^{2h-1} + k$.
\eenum
This completes the proof.  
\end{proof}

 \bt              \label{LinearPerturb:theorem:LowerGrowth-2}
 Let $\varphi = \sum_{i=1}^h c_ix_i$ be a linear form with integer coefficients 
 that satisfies condition $N$.    
There exists an infinite $\varphi$-Sidon set $A$ of  positive integers such that 
\[
A(t) \gg t^{1/(2h-1)}.  
\]
 \et

\begin{proof}
This follows from inequality~\eqref{LinearPerturb:GrowthIneq}.  
\end{proof}

\section{Open problems}             \label{LinearPerturb:section:OpenProblems}
\benum
\item
Let $\varphi = \sum_{i=1}^h c_i x_i$ be a linear form with integer coefficients.  
Let $\PP$ be the set of prime numbers and let $A = \{\log p: p \in \PP\}$.
Consider the $h$-tuple $(p_1,\ldots, p_h) \in \PP^h$ of not necessarily distinct prime numbers, 
and let $\PP_0 = \{p \in \PP : p=p_i \text{ for some } i \in \{1,\ldots, h\} \}$.  
For each $p \in \PP_0$, let 
\[
I_p = \{i\in \{1,\ldots, h\}: p_i = p \} \qqand s(I_p) = \sum_{i\in I_p} c_i.
\]
We have 
\[
\varphi(p_1,\ldots, p_h) = \sum_{i=1}^h c_i \log p_i = \sum_{p\in \PP_0} s(I_p) \log p 
= \log \prod_{p\in \PP_0} p^{S(I_p)}. 
\]
If the linear form $\varphi$ satisfies property $N$, then, by the fundamental theorem of arithmetic, 
the set $A = \{\log p: p \in \PP\}$ is a $\varphi$-Sidon se.  

For the linear form $\psi = x_1 + \cdots + x_h$, Ruzsa~\cite{ruzs98a} used the set $A$ to construct 
large classical Sidon sets of positive integers .  
Are such constructions also possible for $\varphi$-Sidon sets of positive integers? 

\item
Let $A = \{ a_k : k=1,2,3,\ldots \}$ and $B = \{ b_k: k = 1,2,3,\ldots\}$ 
be sequences of integers.  
The set $A$ is a \emph{polynomial perturbation} of  $B$ if 
\[
|a_k - b_k| < k^r 
\] 
for some $r > 0$ and  all $k \geq k_0$.  
The set $A$ is a \emph{bounded perturbation} of  $B$ if 
\[
|a_k - b_k| < m_0
\]
for some $r > 0$ and  all $k \geq k_0$.  

Let $\varphi$ be a linear form with integer coefficients that satisfies condition $N$.  
Let $B$ be a set of integers.  Does there exist a $\varphi$-Sidon set of integers 
that is a polynomial perturbation of $B$?  

Does there exist a $\varphi$-Sidon set of integers 
that is a bounded perturbation of $B$?

\item
Let $\varphi$ be a linear form with integer coefficients that satisfies condition $N$.  
For every positive integer $n$, determine the cardinality of the largest $\varphi$-Sidon 
subset of $\{1,2,\ldots, n\}$.

\item
There exists $c > 0$ such that, for every positive integer $n$, 
 there is a classical Sidon set $A \subseteq \{1,\ldots, n\}$ 
with $A(n) \geq c  \sqrt{n}  $.  
However, there is no infinite classical Sidon set $A$ of positive integers 
such that $A(n) \geq c \sqrt{n} $ for some $c > 0$ and all $n \geq n_0$.  
Indeed, Erd\H os (in St{\"o}hr~\cite{stoh55}) proved that 
every infinite classical Sidon set satisfies 
\[
\liminf_{n\rightarrow \infty} A(n)\sqrt{ \frac{\log n}{n} } \ll 1.
\]
Are there analogous lower bounds for infinite $\varphi$-Sidon sets of positive integers 
associated with binary linear forms $\varphi = c_1x_1 + c_2x_2$ or with linear forms 
$\varphi = \sum_{i=1}^h c_i x_i$ for $h \geq 3$?

\item
Consider  sets of integers.  
One might expect that the elements of a set $A$ of integers that is ``sufficiently random'' or 
``in general position'' will be a  classical Sidon set, that is, will not contain a nontrivial solution 
of the equation $x_1+x_2 = x_3+x_4$.  Equivalently, the set $A$ will be one-to-one (up to transposition) on 
the function $f(x_1,x_2) = x_1 + x_2$.  There is nothing special about the function $x_1 + x_2$.  
One could ask if $A$ is one-to-one  (up to permutation) on some symmetric function, 
or one-to-one on a function that is not symmetric.  
The  functions considered in this paper are linear forms in $h$ variables.  

Conversely, given the set $A$ of integers, we can ask what are the functions 
(in some particular set \mcf\ of functions) with respect to which the set $A$ is one-to-one.  
This inverse problem is considered in Nathanson~\cite{nath2021-202}.

\eenum

\def\cprime{$'$} \def\cprime{$'$}
\providecommand{\bysame}{\leavevmode\hbox to3em{\hrulefill}\thinspace}
\providecommand{\MR}{\relax\ifhmode\unskip\space\fi MR }
\providecommand{\MRhref}[2]{%
  \href{http://www.ams.org/mathscinet-getitem?mr=#1}{#2}
}
\providecommand{\href}[2]{#2}

\end{document}